\documentclass[11pt]{amsart}
\usepackage{amscd} 
\usepackage{amsfonts} 
\usepackage{amssymb} 
\usepackage{latexsym} 

\input diagrams

\newcommand{\ncm}{\newcommand}

\newtheorem{theorem}{Theorem}[section]
\newtheorem{prop}[theorem]{Proposition}
\newtheorem{lemma}[theorem]{Lemma}

\newtheorem{lem&def}[theorem]{Lemma \& Definition}

\newtheorem{example}[theorem]{Example}

\def\C{\mathbb{C}\,}

\ncm{\End}{\mbox{\rm End}\,}

\def\Hom{\mbox{\rm Hom}\,}

\def\id{\mbox{\rm id}}

\def\into{\hookrightarrow}
\def\to{\rightarrow}

\def\o{\otimes}    

\def\ket{\rangle}

\ncm{\rarr}[1]{\stackrel{#1}{\longrightarrow}}
\ncm{\larr}[1]{\stackrel{#1}{\longleftarrow}}
\def\cop{\Delta}

\def\eps{\varepsilon}

\def\du1{\hat 1}

\def\-1{_{(-1)}}
\def\0{_{(0)}}
\def\1{_{(1)}}
\def\2{_{(2)}}
\def\3{_{(3)}}

\def\du1{\hat 1}
\def\lact{\triangleright}
\def\ract{\triangleleft}

\begin{document}

\title[Hopf algebroids and Galois extensions]{Hopf algebroids and Galois extensions}
\author{Lars Kadison} 
\address{University of New Hampshire \\ Kingsbury Hall \\ 
Durham, NH 03824 U.S.A.} 
\email{kadison@math.unh.edu}
\thanks{The author thanks 
NORDAG in Bergen, Eric Grinberg and Andrew Rosenberg at the University of New Hampshire,
and the organizers of the Brussels conference in May 2002 where the author gave
a talk on the subject of sections~2 and~5.}
\subjclass{06A15, 12F10, 13B02, 16W30}  
\footnote{To appear in the \textit{Bulletin of the Belgian Mathematical Society-Simon Stevin}.}
\date{}

\begin{abstract}
To a finite Hopf-Galois extension $A | B$ we associate 
dual
bialgebroids $S := \End\,_BA_B$ and $T := (A \o_B A)^B$ over the centralizer 
$R$  using the depth two theory  
in  \cite[Kadison-Szlach\'anyi]{KSz}. First we extend results on the equivalence of certain
properties of  Hopf-Galois
extensions with corresponding properties of the coacting Hopf algebra  \cite{KT,Doi}
 to depth two extensions using coring theory \cite{Brz}. Next we show 
 that $T^{\rm op}$ is a Hopf algebroid over the centralizer $R$
via Lu's theorem \cite[5.1]{Lu} for smash products with special modules over the
Drinfel'd double, the Miyashita-Ulbrich action, the
fact that $R$ is a commutative algebra in
the pre-braided category of Yetter-Drinfel'd modules \cite{Sch}
and the equivalence of Yetter-Drinfel'd modules with modules over
Drinfel'd double \cite{Maj}. 
In our last section, an exposition of results of Sugano \cite{Su82,Su87}
 leads us  to a Galois correspondence 
between sub-Hopf algebroids
 of $S$ over simple subalgebras of the centralizer with finite projective 
intermediate simple subrings 
of a finite projective H-separable extension of simple rings $A \supseteq B$.
\end{abstract} 
\maketitle

\section{Introduction}

The notion of a Hopf-Galois 
extension  was introduced by Kreimer and Takeuchi in 1981 
\cite{KT} as a generalization of Galois extensions of fields, commutative rings
and noncommutative rings, 
and studied  in connection with
affineness theorems for algebraic groups, non-normal separable field extensions and
Takesaki duality in operator algebras by Schneider,
Greither-Pareigis, Blattner-Montgomery and others. 
 Finite Hopf-Galois extensions 
 have a theory similar to that of depth two finite index 
 subfactors in the von Neumann algebra theory of ``continuous geometry,'' 
the explanation being that both are depth two ring extensions \cite{KN, KSz}. 

Hopf algebroids over noncommutative
rings were introduced by Lu \cite{Lu} in connection with  quantization of Poisson
groupoids in Poisson geometry. Examples of Hopf algebroids are first and foremost
Hopf algebras and groupoid algebras but more significantly come from
solutions to dynamical Yang-Baxter equations \cite{EV}, 
weak Hopf algebras 
\cite{BSz,EN}, finite index subfactors \cite{EV} and in the study of the non-flat
case of index theory for transversally elliptic operators \cite{CM, BSz2}. 

 A bialgebroid $S$, i.e.,  a 
Hopf algebroid without antipode, and its $R$-dual $T$ has been
associated with a depth two ring extension $A | B$ with
centralizer $R$ in Kadison-Szlach\'anyi \cite{KSz}. $S$ 
acts from the left on the over-ring $A$ such that the right
endomorphism ring is isomorphic to a smash product $A \rtimes S$ \cite{KSz}. 
Moreover, $T$ acts from the right 
on the left endomorphism ring $\mathcal{E}$ \cite{KSz} such
 that the endomorphism ring $\End\,_A A\! \o_B\! A$
is similarly isomorphic to a smash product $T \ltimes \mathcal{E}$, which leads
to a Blattner-Montgomery duality result if the extension $A | B$ is
also Frobenius \cite{K2003}. 

In this paper we show via Lu's theorem \cite[5.1]{Lu} 
that the bialgebroid $T^{\rm op}$ of an $H$-Galois extension 
$A$ with subring of invariants
$B$ has  Hopf
algebroid structure over  $R$.  In order to frame it in terms of Lu's hypotheses,
the proof makes use of Miyashita-Ulbrich action, Yetter-Drinfel'd modules
and Drinfel'd doubles.  It is perhaps interesting to mention that Lu's 
theorem  is a quantization of another theorem by Lu in
Poisson geometry \cite[1.2]{Lu90,Lu} via a dictionary between Poisson geometry
and noncommutative algebra \cite{Lu}. In section~3 we establish some theorems
that inform us when depth two extension $A/B$ are separable or Frobenius judging
from the dual properties of the underlying $R$-corings of the acting
bialgebroids $S$ or $T$.  
In a final expository section of this paper, we show 
that a one-sided f.g.\ projective H-separable extension of simple rings,
such as special finite Jones index subfactors with simple relative commutant, enjoys a Galois correspondence between
intermediate simple rings forming f.g. projectives with the overring,
and Hopf subalgebroids over the simple subalgebras of the centralizer.  
This  depends on Sugano's one-to-one correspondence
between the intermediate simple subrings and simple subalgebras of the
centralizer of the full H-separable extension \cite{Su82,Su87}, with its roots in
work on certain classical
inner Galois theories of simple artinian rings and division rings
 by Jacobson, Bourbaki, Tominaga and others.     We hope that this exposition
will be a first step toward an algebraic generalization of the Galois correspondence by
Nikshych and Vainerman between finite depth and index intermediate subfactors
and coideal subalgebras of a weak $*$-Hopf algebra \cite{NV}.


\section{Dual bialgebroids over the centralizer}  

In this section we review the basics of the dual bialgebroid constructions in \cite{KSz},
while computing the bialgebroids of a finite Hopf-Galois extension as
a running example.  

Let $B$ be a unital subring of $A$, an associative noncommutative ring with unit,
or an image of a ring homomorphism $\overline{B} \to A$.   Recall
that the ring extension $A | B$ is said to be of \textit{depth two} if
$$A \o_B A \oplus * \cong \oplus^n A$$ as natural $B$-$A$ and $A$-$B$-bimodules
\cite{KSz}.  Equivalently, there are elements $\beta_i \in S := \End\,_BA_B$,
$t_i \in T := (A \o_B A)^B$ (called a \textit{left D2 quasibasis}) such that ($a,a' \in A$)
\begin{equation}
\label{eq: left d2 qb}
 a \o a' = \sum_{i=1}^n t_i\beta_i(a)a', 
\end{equation}
and a right D2 quasibasis $\gamma_j \in S$, $u_j \in T$ such that
\begin{equation}
\label{eq: right d2 qb}
 a \o a' = \sum_j a \gamma_j(a') u_j. 
\end{equation}
Fix both D2 quasibases
for our work in this paper.

\begin{example}
\label{ex-rd2}
\begin{sl}
Consider a Hopf-Galois extension $A | B$ with
$n$-dimensional Hopf $k$-algebra $H$ \cite{KT} with
$k$ an arbitrary field.  Our convention is that $H^*$ acts from the 
left on $A$ with subalgebra of
invariants $B$, or equivalently, there is 
a dual right coaction $A \to A \o_k H$, $a \mapsto a\0 \o a\1$:
 the Galois isomorphism $\beta: A \o_B A \stackrel{\cong}{\longrightarrow}
A \o_k H$ given by $\beta(a \o a') = a{a'}\0 \o {a'}\1$ , 
which is an $A$-$B$-bimodule, right $H$-comodule morphism.
It follows that $A \o_B A \cong \oplus^n A$ as $A$-$B$-bimodules.  As 
$B$-$A$-bimodules there is  similarly an isomorphism $A \o_B A \cong \oplus^n A$ 
 by making use of the opposite 
 Galois isomorphism $\beta'$ given by  $\beta'(a \o a') = a\0 a' \o a\1$.

Now compute a right D2 quasibasis $\{ \gamma_i \}$,
$\{ u_i \}$  for $A | B$.  Let $\{ h_i \}$, $\{ p_i \}$
be dual $k$-bases in $H$, $H^*$, respectively. Define $\gamma_i \in \End\,_BA_B$
by $\gamma_i(a) := p_i \cdot a$ ($a \in A$, $i \in \{ 1,2,\ldots,n \}$).
Let $u_i := \beta^{-1}(1 \o h_i) \in (A \o_B A)^B$.  We verify this: ($a,a' \in A$)
\begin{eqnarray*}
\sum_i a \gamma_i(a')u_i & = & \sum_i a (p_i \cdot a')\beta^{-1}(1 \o h_i) \\
& = & \sum_i a {a'}\0  p_i( {a'}\1) \ket  \beta^{-1}(1 \o h_i) \\
& = &  \beta^{-1}(a{a'}\0 \o {a'}\1)  \\
& = & a \o a'
\end{eqnarray*}
\end{sl}
\end{example}

The paper  \cite{KSz} found a bialgebroid with action and smash product structure 
within  the Jones construction above a depth two ring extension $A | B$. Namely, 
if $R$ denotes the centralizer of $B$ in $A$,
a left $R$-bialgebroid structure on $S$ is given by the composition
ring structure on $S$ with source and target mappings
corresponding to the left regular representation $\lambda: R \to S$
and right regular representation $\rho: R^{\rm op} \to S$, respectively. Since these
 commute
($\lambda(r) \rho(r') = \rho(r')\lambda(r)$ for every $r,r' \in R$),
we may induce an $R$-bimodule structure on $S$ solely
  from the left  
by $$r \cdot \alpha \cdot r' := \lambda(r)\rho(r') \alpha = r\alpha(?)r'.$$

Now an $R$-coring structure $\mathcal{S} = (S, \cop, \eps)$ is given
by 
\begin{equation}
\label{eq: cop A}
\cop(\alpha) := \sum_i \alpha(\, ?\, t_i^1) t_i^2 \o_R \beta_i 
\end{equation}
for every $\alpha \in S$, denoting $t_i = t_i^1 \o t_i^2 \in B$ by
suppressing a possible summation, and
\begin{equation}
\eps(\alpha) = \alpha(1)
\end{equation}  
satisfying the additional axioms of a bialgebroid (cf.\ section~4),
such as multiplicativity of $\cop$ and a condition that makes sense of this
requirement. We have the equivalent formula for the coproduct \cite[Th'm 4.1]{KSz}:
\begin{equation}
\label{eq: right D2 cop}
\cop(\alpha) := \sum_j \gamma_j \o_R u_j^1 \alpha(u_j^2\, ?\, ) 
\end{equation}
Since $S \o_R S \cong \Hom_{B-B}(A \o_B A, A)$ via $\alpha \o \beta \mapsto
(a \o a' \mapsto \alpha(a)\beta(a'))$, we have the simpler formula
via identication,
\begin{equation}
\label{eq: Lu1}
\cop(\alpha)(a \o a') = \alpha(aa'),
\end{equation}
which clearly shows this bialgebroid structure on $S$ to be a generalization
to depth two ring extensions of Lu's bialgebroid ${\rm End}_k\, C$ over
a finite dimensional $k$-algebra $C$ (cf.\ section~4).

\begin{example}
\label{ex: ess}
\begin{sl}
We determine the $R$-bialgebroid
$S$ for the Hopf-Galois extension $A | B$ introduced above.
It is well-known (see for example \cite{Mont}) that the right endomorphism
ring is a smash product: 
\begin{equation}
A \rtimes H^* \cong \End A_B
\end{equation}
via $a \rtimes p \mapsto \lambda(a) \circ (p \cdot \, ? \, )$.
This is an $A$-$B$-isomorphism (where $a'(a \o p)b := a'ab \o p$
and $\End A_B$ is the natural $A$-$A$-bimodule).  The $B$-centralizer in $\End A_B$
is of course $(\End A_B)^B = S$, whence
\begin{equation}
\Phi: S \stackrel{\cong}{\longrightarrow} R \rtimes H^*
\end{equation}
with multiplication given by the smash product:
\begin{equation}
\label{eq: smash}
 (r \rtimes p)(r' \rtimes p') = r ( p\1 \cdot r') \rtimes p\2 p'.
\end{equation}
If $t \in H, T \in H^*$ denote a dual pair of left integrals (where
 $T(t) = 1$), and $\sum_i x_i \otimes y_i
= \beta^{-1}(1 \o t)$, $\Phi(\alpha) = \sum_i (\alpha(x_i) \rtimes T)(y_i \rtimes 1)$
for $\alpha \in S$ (cf.\ \cite{Mont}).  

The induced $R$-coring structure is (the trivial structure except for
the more complex right $R$-module action)  given by 
$\tilde{s}(r) = r \rtimes 1$, 
$$\tilde{t}(r) = \Phi(\rho(r)) = \sum_i x_i r (T\1 \cdot y_i) \rtimes T\2, $$
with coproduct 
\begin{eqnarray}
\cop(r \rtimes p) & = & (\Phi \o \Phi)(\sum_i\gamma_i(?)\o_R u_i^1 r(p\cdot(u_i^2 \, ? \, ))
\nonumber \\
& = & (\Phi \o \Phi)(\sum_i \gamma_i(?)u_i^1 r (p\1 \cdot u_i^2)\o_R (p\2 \cdot \, ? \, ) \\
& = &  r \rtimes p\1 \otimes p\2,
\end{eqnarray}
and counit $$\eps(r \rtimes p) = r(p \cdot 1_A) = r \eps(p).$$
The formula for $\cop$ makes use of the depth two eq.~(\ref{eq: right d2 qb}). 
\end{sl}
\end{example}

The left action of $S$ on $A$ is very simply given by evaluation, 
\begin{equation}
\label{eq: alpha-action}
\alpha \lact a = \alpha(a).
\end{equation}
This action has invariant subring (of elements
$a \in A$ such that $\alpha\lact a = \eps(\alpha)a$) equal precisely to $B$ 
if the natural module $A_B$ is balanced \cite{KSz}. 
This action is measuring because $\alpha\1(a)\alpha\2(a') = \alpha(aa')$ by
eq.~(\ref{eq: Lu1}). 

The smash product $A \rtimes' S$, which is $A \o_R S$ as  abelian groups with associative
multiplication given by eq.~(\ref{eq: smash}),
is isomorphic
as rings to $\End A_B$ via $a \o_R \alpha \mapsto \lambda_a \alpha$ \cite{KSz}.

\begin{example}
\begin{sl}
For the $H$-Galois extension $A | B$ just introduced, the action of $S$ on $A$
under the isomorphism $S \cong R \rtimes H^*$ is just given by $(r \rtimes p)\cdot a
= r(p \cdot a)$.  The smash product of $A$ with the bialgebroid $R \rtimes H^*$
just recovers the ordinary smash product of $A$ with $H^*$:
$$ A \rtimes' ( R \rtimes H^*) \cong A \rtimes H^*$$
as ring isomorphism by an easy exercise.  
\end{sl}
\end{example}

For any subring $B$ in ring $A$, the construct
 $T = (A \o_B A)^B$ (``the $B$-central tensor-square of $A$ over $B$'')
 has a unital ring structure induced from 
$T \cong \End\,_A(A\! \o_B\! A)_A  $ via $F \mapsto F(1 \o 1)$, which is given
by 
\begin{equation}
\label{eq: T-mult}
tt' = {t'}^1 t^1 \o t^2 {t'}^2
\end{equation} 
for each $t,t' \in T$. There are
obvious commuting homomorphisms of $R$ and $R^{\rm op}$ into $T$ given by 
$r \mapsto 1 \o r$ and $r' \mapsto r' \o 1$, respectively. From the right,
these two ``source'' and ``target'' mappings induce the $R$-$R$-bimodule
structure ${}_RT_R$ given by $$r \cdot t \cdot r' =  (t^1 \o t^2)(r \o r') =
rtr',$$
the ordinary bimodule structure on a tensor product.  

For a D2 extension $A/B$, 
there is a right $R$-bialgebroid structure on $T$ with coring structure
$\mathcal{T} = (T, \cop, \eps)$ given by the two equivalent formulas: 
\begin{equation}
\label{eq: copB}
\cop(t) = \sum_i t_i \o_R (\beta_i(t^1) \o_B t^2) = \sum_j (t^1 \o_B \gamma_j(t^2)) \o_R
u_j 
\end{equation}
\begin{equation}
\label{eq: epsB}
\eps(t) = t^1 t^2
\end{equation}
By \cite[Th'm 5.2]{KSz} $\cop$ is multiplicative and the other axioms of
a right bialgebroid are satisfied. Since the D2 conditions yield
 $T \o_R T \cong (A \o_B A \o_B A)^B$,
the coproduct enjoys a Lu generalized formula,
\begin{equation}
\label{eq: Lu2}
\cop(t) = t^1 \o 1 \o t^2 \ \ \ (t \in T).
\end{equation} 
Indeed, $T$ is a right-handed generalization of Lu's bialgebroid $C^e = C \o_k C^{\rm op}$
for a finite dimensional $k$-algebra $C$, although $T$, unlike $C^e$, 
has in general no antipode.

\begin{example}
\label{ex: tee}
\begin{sl}
We return to the example of $A | B$
an $H$-Galois extension, to compute the $R$-bialgebroid $T$.
Since $\beta: A \o_B A \to A \o H$ is an $A$-$B$-bimodule isomorphism,
it follows that $T = A^B \o H \cong R \o H$ via $\beta$. We next study
the multiplication $\star$ imposed on $R \o H$ by $\beta$
and the multiplication (\ref{eq: T-mult}) on $T$.  
Let $h,h' \in H$ and $t,t' \in T$ such that $\beta(t) = 1 \rtimes h$
and $\beta(t') = 1 \rtimes h'$.  
We compute using the fact that $\beta$ is an $H$-comodule homomorphism in
the last step:  ($r,r' \in R$, $h,h' \in H$) 
\begin{eqnarray}
(r \o h) \star(r' \o h') & = & \beta(rt) \beta(r't') \nonumber \\
& = & \beta(r' {t'}^1 r t^1 \o t^2 {t'}^2) \nonumber \\
& = & r' {t'}^1 r (t^1 {t^2}\0){{t'}^2}\0 \o {t^2}\1 {{t'}^2}\1 \nonumber \\
& = & r' {b'}^1 r {{b'}^2}\0 \o h {b'}\2 \nonumber \\
& = & r'(r \ract {h'}\1) \o h {h'}\2
\label{eq: mu-mult}
\end{eqnarray}
where $\ract$ denotes the Miyashita-Ulbrich action of $H$ on $R$ from
the right \cite{U,DT,Sch,K2001}.  (Recall that if $\beta(t) = 1 \o h$
then $r \ract h := t^1 r t^2$.) From this formula for $\star$, we
see that $\beta$ induces an algebra isomorphism, 
\begin{equation}
T^{\rm op} \cong R \rtimes H^{\rm op}
\end{equation}
where the right action by $H$ is equivalent to a left action by $H^{\rm op}$.

The $R$-coring structure on $R \rtimes H^{\rm op}$ induced from
$T^{\rm op}$ is (the trivial structure) given by ($b := \beta^{-1}(1 \o h)$)
\begin{eqnarray}
\tilde{s}(r) & = & r \rtimes 1 \ \ (r \in R) \\
\tilde{t}(r) & = & r\0 \rtimes r\1  \\
\cop( r \rtimes h) & = & r \beta \o \beta \cop_T (\beta^{-1}(1 \o h) \nonumber \\
                  & = & r \beta(b^1 \o {b^2}\0 p_i({b^2}\1)  \o \beta(u_i) \nonumber \\ 
                  & = & r \beta(b^1 \o {b^2}\0) \o {b^2}\1  \nonumber \\
                  & = & r \rtimes h\1 \o h\2 \\
\eps(r \rtimes h) & = & \eps_T(\beta^{-1}(r \o h) \nonumber \\
              & = & r \eps(h)
\end{eqnarray}
The formula for $\cop$ again uses the right $H$-comodule property of $\beta$,
while the formula uses the counitality of the $H$-comodule $A$ with Eq.\ (\ref{eq: epsB}). 
\end{sl}
\end{example}

\begin{example}
\begin{sl}
The first of two special cases of finite Hopf-Galois extensions with normal
basis property is naturally a finite
dimensional Hopf algebra $H$ coacting on itself via its comultiplication $\cop$. 
The coinvariant subalgebra $B$ is the unit subalgebra $k1_H$, $R = H$, the $R$-bialgebroid 
$S$ is ${\rm End}_k\, H \cong H \rtimes H^*$ by example~\ref{ex: ess}, and the
bialgebroid structure is the same as  the ``Heisenberg double'' of $H$ in Lu's
 \cite[section~6]{Lu},  for which Lu finds an antipode and Hopf algebroid structure. 

The Miyashita-Ulbrich action of $H$ on itself from the right is given by
ordinary conjugation, $h \ract a = S(a\1) h a\2$.  Thus the $H$-bialgebroid structure
on $T^{\rm op}$ is given above in example~\ref{ex: tee} --- with antipode and
Hopf algebroid structure in section~5 below.  

The second example of an elementary nature is obtained from groups $G$ and $N$ where
$N$ is a normal subgroup of $G$ of index $n$ and $G/N$ its factor group. Given any field $k$, the group algebra $A = k[G]$
is Galois over $B = k[N]$ with cocommutative Hopf algebra $H = k[G/N]$. The Galois
map $\beta: A \o_B A \to A \o   H$ is given by $\beta(g \o g') = gg' \o g'N$ for every
$g,g' \in G$. Given a set of right coset representatives $g_1,\ldots,g_n$,
the prescription for finding  right D2
quasibases
 in example~\ref{ex-rd2} yields  $u_i = g_i^{-1} \o g_i$ 
and $\gamma_i(g) = 0$ if $gN \neq g_iN$
and $\gamma_i(g) = g$ if $gN = g_iN$.  Since $\beta^{-1}(1 \o gN) = g^{-1} \o g$,
the action associated to $T^{\rm op}$ above is the Miyashita action  given by 
$x \ract gN = g^{-1} x g$ where $x \in C_A(B)$.  
\end{sl}
\end{example}

\section{When D2 extensions are separable, split or Frobenius}

Given a D2 extension $A / B$, we made the acquaintance in the previous section
 of the underlying $R$-corings $\mathcal{S}$
and $\mathcal{T}$ of the $R$-bialgebroids $S = {\rm End}_{B-B}\, A$ and $T
= (A \o_B A)^B$, respectively.  In this section we show that coring properties of $\mathcal{S}$
or $\mathcal{T}$  such as coseparability determine properties of $A/B$ such
as separability, and vice versa.

For the next theorem, recall that any $R'$-coring $(\mathcal{C}, \cop, \eps)$ is {\it coseparable}
if there is an $R'$-$R'$-homomorphism $\gamma: \mathcal{C} \o_{R'} \mathcal{C} \to R'$
(called a {\it cointegral}) 
such that $\gamma(c\1 \o c\2) = \eps(c)$ and $c\1 \gamma(c\2 \o c') = \gamma(c \o c'\1)c'\2$
for every $c,c' \in \mathcal{C}$ (cf.\ \cite{Brz,BKW,CMZ}).

\begin{theorem}
Let $A/B$ be a right f.g.\ projective D2 extension.  Then $A / B$ is a separable extension
if and only if the $R$-coring $\mathcal{S}$ is coseparable.
\end{theorem}
\begin{proof}  
($\Rightarrow$ \cite[Example 3.6]{BKW}) Given separability element $e = e^1 \o e^2 \in (A \o_B A)^A$ for $A/B$,
define cointegral $\gamma: \mathcal{S} \o_R \mathcal{S} \to R$ by
$\gamma(\alpha \o \beta) = \alpha(e^1)\beta(e^2)$.  The rest of the proof 
follows \cite[Example 3.6]{BKW} and does not require $A_B$ to be finite projective.

Suppose a dual basis for the natural module $A_B$ is given by $\{ a_k \}$,
$\{ f_k \} $.  

($\Leftarrow$)  Given cointegral $\gamma: \mathcal{S} \o \mathcal{S}  \to R$, define
$e = \sum_i b_i \gamma(\beta_i \o I_A)$ where $I_A$ is the identity map on $A$
and $\{ b_i \}$, $\{ \beta_i \}$ is the left D2 quasibases introduced above.
Of course, $e \in (A \o_B A)^B$; also,  ($\alpha \in \mathcal{S}$)
$$ \alpha(e^1) e^2 =  \sum_i \alpha(b^1_i) b^2_i \gamma(\beta_i \o I_A) = \gamma(\sum_i 
\alpha(b_i^1) b_i^2\beta_i \o I_A)
= \gamma(\alpha \o I_A),$$
whence if $\alpha = I_A$, $e^1 e^2 = \gamma(I_A \o I_A)  = \eps(I_A) = 1_A$  
since $\cop(I_A) = I_A \o I_A$.  

It follows that $\alpha(e^1)e^2a = \gamma(\alpha \o I_A)a$ for $a \in A$, but
$$\alpha(ae^1)e^2 =  \alpha\1(a) \gamma(\alpha\2 \o I_A)= \gamma(\alpha \o I_A)I_A(a)
= \alpha(e^1)e^2a.$$
Since $ A \rtimes S \cong \End A_B $ via $a \rtimes \alpha \mapsto \lambda(a) \circ 
\alpha$, it follows that $f(e^1)e^2a = f(ae^1)e^2$ for each $f \in \End A_B$. Finally then computing in $A \o_B A$:
$$ae = \sum_k a_k \o f_k(ae^1)e^2 = \sum_k a_k f_k(e^1) \o e^2a = ea, $$
for each $a \in A$, whence $e$ is a separability element of $A/B$. 
\end{proof}

\begin{example}
\begin{sl}
Suppose again that $A/B$ is an $H$-Galois extension.  The multiplication
mapping $A \o_B A \to A$ corresponds under the Galois isomorphism $\beta$
to $A \o_k H \to A$ given by $a \o h \mapsto a\eps(h)$.
It follows from the theorem that $A/B$ is a separable extension iff $H$ is
semisimple, since $H$ is semisimple iff the left integral $t \in H$
may be chosen so that $\eps(t) = 1$, whence $e = \beta^{-1}(1 \o t)$
is a separability element.  This recovers a theorem of Doi \cite{Doi}.
\end{sl}
\end{example}

For the next theorem, we recall that any $R$-coring $\mathcal{C}$ is {\it cosplit}
if there is $e \in \mathcal{C}^{R'}$ such that $\eps(e) = 1$, i.e., the counit
$\eps: \mathcal{C} \to R'$ is a split $R'$-$R'$-epi. An ring extension
$A'/B'$ is {\it split} if there is a $B'$-$B'$-epimorphism $E: A' \to B'$ such that
$E(1) = 1$ (cf.\ \cite{Brz,CMZ}). 

\begin{example}
\begin{sl}
If $A/B$ is a split extension, the Sweedler $A$-coring $A \o_B A$ \cite{Sw}
is coseparable \cite{Brz}; similarly one shows that if $A/B$ is D2 and split, $\mathcal{T}$
is coseparable.

If $A/B$ is separable and D2, then $\mathcal{T}$ is a cosplit $R$-coring, since
a separability element $e \in \mathcal{T}$ satisfies $\eps(e) = e^1e^2 = 1$
and $e \in \mathcal{T}^R$. Define a ring extension $A/B$ to be {\it Procesi}
if $BR = A$; e.g., centrally projective extensions or extensions of commutative
rings are Procesi. Conversely then, $\mathcal{T}$ cosplit implies
$A/B$ is separable if $A/B$ is a D2 Procesi ring extension. 
\end{sl}
\end{example}

\begin{theorem}
Suppose $A/B$ is a D2 extension with double centralizer condition $C_A(C_A(B)) = B$.
Then  $A/B$ is a split extension iff $\mathcal{S}$ is a cosplit $R$-coring.
\end{theorem}
\begin{proof}  The proof only requires $C_A(R) = B$ in the direction $\Leftarrow$. 

($\Rightarrow$)  If $E: A \to B$ splits the inclusion map, then $E \in \mathcal{S}^R$
since $rE(a) = rE(a)$ for each $r \in C_A(B), a \in A$.  Moreover, $\eps(E) = E(1) = 1$
and we conclude $\mathcal{S}$ is cosplit.

($\Leftarrow$) Suppose $e \in \mathcal{S}^R$ such that $e(1) = 1$. Since $e(a)r = re(a)$
for $a \in A$, $e(a) \in C_A(R) = B$, whence $e: A \to B$ splits the inclusion $B \into A$.
\end{proof}

\begin{example}
\begin{sl}
As noted in \cite{K2003}, an H-separable extension is D2.  
If $A/B$ is an H-separable extension and $A_B$ is balanced, then 
$A/B$ is  D2  and $C_A(C_A(B)) = B$: see Lemma~\ref{lemma-seeme}.

Another example:  if $A/B$ is an H-separable extension of simple rings with
$A_B$ f.g. projective, then $A/B$ is D2 and $C_A(C_A(B)) = B$.  
(Cf.\ Prop.~\ref{prop: sugano}.)  

It is a problem which would generalize and improve 
 results of Noether-Brauer-Artin on simple rings, 
if $A/B$ a right progenerator H-separable
extension implies $A/B$ is split \cite{Su1999}.
\end{sl}
\end{example}

Recall that an $R'$-coring $\mathcal{C}$ is {\it Frobenius} if
there is an $R'$-$R'$-coring $\gamma: \mathcal{C} \o_{R'} \mathcal{C}
\to R'$ and $e \in \mathcal{C}^{R'}$ such that $\gamma(c \o e) = \eps(c) =
\gamma(e \o c)$ and $\gamma(c \o c'\1)c'\2 = c\1 \gamma(c\2 \o c')$ for every
$c,c' \in \mathcal{C}$ (cf.\ \cite{Brz,CMZ}).  

\begin{prop}
Let $A/B$ be a D2 right progenerator Procesi extension.  Then $A/B$ is a Frobenius extension iff
$\mathcal{T}$ is a Frobenius coring.
\end{prop} 
\begin{proof}
($\Rightarrow$) Suppose $(E: A \to B, e \in (A \o_B A)^A)$ is a Frobenius system
for $A/B$; i.e., for each $a \in A$, we have $a = e^1 E(e^2 a) = E(ae^1)e^2$.
Define $\gamma: \mathcal{T} \o_R \mathcal{T} \cong (A \o_B A \o_B A)^B \to R$
by $\gamma(a \o a' \o a'') = a E(a') a'' \in R$.
It follows that: ($d \in \mathcal{T}$)
$$ \gamma(d \o e) = \gamma(d^1 \o d^2 e^1 \o e^2) = d^1 E(d^2 e^1)e^2 = d^1 d^2 = \eps(d),$$
and similarly $\gamma(e \o d) = \eps(d)$.
Recalling the $E$-multiplication $\cdot_E$ on $A \o_B A$ induced by 
$A \o_B A \cong \End A_B$, we note:
\begin{eqnarray*}
d\1 \gamma(d\2 \o d') & = & \sum_i b_i \beta_i(d^1) E(d^2 d'^1) d'^2 \\
                      & = & d^1 \o_B E(d^2 {d'}^1) {d'}^2 = d \cdot_E d' \\
                      & = & \sum_j \gamma(d \o_R d'^1 \o_B \gamma_j(d'^2))c_j \\
                      & = & \gamma(d \o d'\1) d'\2
\end{eqnarray*}

($\Leftarrow$) 
We now assume that $BR = A$ and that $A_B$ is a progenerator. 
We see from \cite[3.3.10]{CMZ}  that $R \to \mathcal{T}^*$
given by $r \mapsto r\eps$ is a Frobenius extension.
But the $R$-dual $\mathcal{T}^* \cong \mathcal{S}$ via $\phi \mapsto \sum_i \phi(b_i)\beta_i$ with inverse
$$ \alpha \mapsto (t \mapsto \alpha(t^1)t^2). $$ The composite $R \to \mathcal{S}$ is the left regular
map $\lambda: R \to \mathcal{S}$, which is therefore Frobenius. Let $\rho_k$, $\eta_k \in \mathcal{S}$,
$E': \mathcal{S} \to R$ be a Frobenius system satisfying ($\alpha \in \mathcal{S}, r,r' \in R$)

\begin{eqnarray}
\sum_k \lambda(E'(\alpha \rho_k)) \eta_k & = \alpha & = \sum_k \rho_k \lambda(E'(\eta_k \alpha)) \label{eq: uno} \\
E'(\lambda(r) \alpha \lambda(r')) & = & rE'(\alpha)r'
\end{eqnarray}
The last equation is equivalent to $E'(\lambda(r \alpha\1(r')) \alpha\2) = r \alpha\1(r') E'(\alpha\2) = rE'(\alpha)r'$.
Since $\alpha(r)b = b\alpha(r)$ for all $b \in B, r \in R, \alpha \in \mathcal{S}$, it follows that for
every $a, a' = \sum_i b_i r_i \in A$
\begin{equation}
\label{eq: vi}
a\alpha\1(a')E'(\alpha\2) = \sum_i a \alpha\1(r_i) E'(\alpha\2)b_i = aE'(\alpha)a' 
\end{equation}
where of course $b_i \in B, r_i \in R$.  

We now claim $\lambda: A \into \End A_B$ is a Frobenius extension, from which it follows that $A/B$ is Frobenius
by a converse of the endomorphism ring theorem, since $A_B$ is a progenerator
 \cite{K99}. This follows from $A \rtimes \mathcal{S} \cong \End A_B$
via $a \rtimes \alpha \mapsto \lambda(a) \alpha$ 
and the assumption $A = BR$.  Define   $E: \End A_B \to A$ by $E(\lambda(a)\alpha) = a E'(\alpha)$.  
Now $E \in {\rm Hom}_{A-A}(\End A_B, A)$ by eq.~(\ref{eq: vi}), since $\lambda(a) \alpha \lambda(a')=
\lambda(a \alpha\1(a')) \alpha\2$. It follows that
 $$\sum_k E(\lambda(a) \alpha \rho_k)\eta_k = \lambda(a) \alpha $$
by eqs.~(\ref{eq: uno}), which also imply (assuming $a = \sum_j b_j r_j \in BR$)
$$ \sum_k \rho_k E(\eta_k \lambda(a) \alpha) =
\sum_{k,j} \lambda(r_j) \rho_k \lambda(b_j) E'(\eta_k \alpha) =
\lambda(a) \alpha. \qed $$
\renewcommand{\qed}{}\end{proof}


\section{Hopf algebroids}

 For the convenience of the reader and the sake of convention, let's recall some facts about
Lu's Hopf algebroid, which consists of a  
left bialgebroid $(H,R,\tilde{s},\tilde{t},\cop, \eps)$,
and and  
\textit{antipode} $\tau$ for $H$.  $H$ and $R$ are $k$-algebras and all 
maps are $k$-linear.  First,
recall from \cite{Lu} (and compare \cite{BM, KSz}) that
the \textit{source} and \textit{target} maps $\tilde{s}$ and $\tilde{t}$ 
are algebra homomorphism and anti-homomorphism, respectively,
of $R$ into $H$ such that $\tilde{s}(r)\tilde{t}(r') = \tilde{t}(r')\tilde{s}(r)$ 
for all $r,r' \in R$.
This induces an $R$-$R$-bimodule structure on $H$ (from the left in this case) by
$r \cdot h \cdot r' = \tilde{s}(r)\tilde{t}(r') h$ ($h \in H$).  With respect to this bimodule
structure, $(H,\cop,\eps)$ is an $R$-coring (cf.\ \cite{Sw}), i.e. with coassociative
coproduct and $R$-$R$-bimodule map $\cop: H \to H \o_R H$ and
counit $\eps: H \to R$ (also an $R$-bimodule mapping).  The image
of $\cop$, written in Sweedler notation, is required to satisfy
\begin{equation}
 a\1 \tilde{t}(r) \o a\2 = a\1 \o a\2 \tilde{s}(r)
\end{equation}
for all $a \in H, r \in R$.  It then makes sense to require that $\cop$ be homomorphic:
\begin{equation}
\cop(ab) = \cop(a) \cop(b), \ \ \ \cop(1) = 1 \o 1
\end{equation}
for all $a,b \in H$. The counit must satisfy the following modified augmentation law:
\begin{equation}
\eps(ab) = \eps(a s(\eps(b))) = \eps(a t(\eps(b))), \ \ \eps(1_H) = 1_R. 
\end{equation}

The axioms of a right bialgebroid $H'$ are opposite those of a left bialgebroid in
the sense that $H'$ obtains its $R$-bimodule structure from the right via its
source and target maps and, from the left bialgebroid $H$ above, we have that 
 $(H^{\rm op}, R, \tilde{t}^{\rm op}, \tilde{s}^{\rm op}, \cop, \eps)$
in this  precise order is a right bialgebroid: 
for the explicit axioms, see \cite[Section~2]{KSz}. 

The left $R$-bialgebroid $H$
is a \textit{Hopf algebroid}  $(H,R,\tau)$ 
if  $\tau: H \to H$ is an algebra anti-automorphism (called an {\it antipode}) such that 
\begin{enumerate}
\item $\tau \tilde{t} = \tilde{s}$;
\item $\tau(a\1)a\2 = \tilde{t} (\eps (\tau(a)))$ for every $a \in A$;
\item there is a linear section $\eta: H \o_R H \to H \o_K H$ to the natural
projection $H \o_k H \to H \o_R H$ such that:
$$ \mu (H \o \tau) \eta \cop = \tilde{s} \eps. $$
\end{enumerate}
The following lemma covers some examples in the literature (e.g. \cite[3.2]{Khal}). 

\begin{lemma}
If  $(H,R,s,t,\cop,\eps,\tau)$ and $(H',R',s',t',\cop',\eps',\tau')$ are Hopf algebroids,
then $$(H \o H',\, R \o R',\, s\o s',\, t \o t',\, (1 \o \sigma \o 1)\cop \o \cop',\,
 \eps \o \eps',\, \tau \o \tau')$$ is (the tensor) Hopf algebroid.  
\end{lemma}
\begin{proof}
The proof is straightforward and left to the reader, $\sigma$ denoting the twist and 
the linear section being given up to two twists by
$\eta \o \eta'$ if $\eta,\eta'$ are the sections for $H$ and $H'$ as in
axiom (3) above. 
\end{proof}

Lu's examples of bialgebroids and Hopf algebroids are the 
following.  Given an algebra $C$ over commutative ground ring $K$ such that
$C$ is finitely generated projective as $K$-module, the following two are left
bialgebroids over $C$ (with $\o = \o_K$):

\begin{example}
\label{ex: Lu BA}
\begin{sl}
 The endomorphism algebra $E := {\rm End}_K\, C$ with $\tilde{s}(c) = \lambda(c)$, $\tilde{t}(c') = 
\rho(c')$, coproduct $\cop(f)(c \o c') = f(cc')$ for $f \in \End_K C$
after noting that $E \o_C E \cong {\rm Hom}_K\, (C \o C, C)$ via $f \o g \mapsto
(c \o c' \mapsto f(c)g(c'))$.  The counit
is given by $\eps(f) = f(1)$.  We see that this is the left bialgebroid
$S$ above when $B = K$, a subring in the center of $A =  C$. 
\end{sl}
\end{example} 

\begin{example}
\label{ex: Lu HA}
\begin{sl}
 The ordinary tensor algebra
$C \o C^{\rm op}$ with $\tilde{s}(c) = c \o 1$, $\tilde{t}(c') = 1 \o c'$ with
bimodule structure $c \cdot c' \o c'' \cdot c''' = cc' \o c'' c'''$.
Coproduct  $\cop(c \o c') = c \o 1 \o c'$ after a simple identification,
with counit $\eps(c \o c') = cc'$ for $c,c' \in C$. $C \o C^{\rm op}$ is
a left $C$-bialgebroid by arguing as in \cite{Lu}, or \cite[$N = K$]{KSz} since
$C | K$ is D2. In addition, $\tau: C \o C^{\rm op} \to C \o C^{\rm op}$ defined
as the twist $\tau(c \o c') = c' \o c$ is an antipode satisfying the axioms
of a Hopf algebroid (in addition,  $\tau^2 = \id$, an \textit{involutive}
antipode). 
\end{sl} 
\end{example}
 
A bialgebroid homomorphism from $(H_1,R_1,s_1, t_1,\cop_1,\eps_1)$
into $(H_2,R_2,s_2$, \newline
$t_2,\cop_2,\eps_2)$ consists of a pair of algebra homomorphisms,
$F: H_1 \to H_2$ and $f: R_1 \to R_2$, such that four squares commute: 
$Fs_1 = s_2 f$,
$Ft_1 = t_2 f$, $\cop_2 F = p(F \o F)\cop_1$ and $\eps_2 F = f \eps_1$, 
where  $f$ induces an $R_1$-$R_1$-bimodule structure on $H_2$ via 
``restriction of scalars,''  $p: H_2 \o_{R_1} H_2 \to H_2 \o_{R_2} H_2$
is the canonical mapping and $F: {}_{R_1}{H_1}_{R_1} \to {}_{R_1}{H_2}_{R_1}$
is a bimodule homomorphism since
$$F(r \cdot h \cdot r') = F(s_1(r) t_1(r') h) = s_2(f(r))t_2(f(r'))F(h) = 
r \cdot_f F(h) \cdot_f r'. $$

If $F$ and $f$ are both inclusions, we say $H_1$ is a \textit{sub-bialgebroid}
of $H_2$; if moreover $H_1$ and $H_2$ are both Hopf algebroids with antipodes $\tau_1$
and $\tau_2$ such that $F \tau_1 = \tau_2 F$, we call $H_1$ a \textit{Hopf subalgebroid}
of $H_2$.  We say that a bialgebroid $H$ is minimal over its base ring $R$ if
it has no proper $R$-subbialgebroid.  
 

\section{$T^{\rm op}$ is a Hopf algebroid}
\label{sec: Hopf algebroids}

In this section, we find an antipode for the bialgebroid $T^{\rm op}$
 we associated to
the $H$-Galois extension $A | B$ in Section~2.  We apply \cite[Theorem 5.1]{Lu}, 
repeated below without proof for the convenience of the reader, 
after noting that the centralizer $R = C_A(B)$ is a commutative algebra in the 
Yetter-Drinfel'd category $\mathcal{YD}_H^H$ of modules-comodules over $H$
\cite[3.1]{Sch}. 

\begin{theorem}[ Lu Theorem 5.1 \cite{Lu}]
\label{thm-Lu}
Let $H'$ be a Hopf algebra with antipode $\tilde{S}$ and $D(H')$ its Drinfel'd double.  Let $V$ be
a left $D(H')$-module algebra.  Assume that the $R$-matrix  $\sum_i (1 \o h_i) \o 
(p_i \o 1)$ satisfies the following pre-braided commutativity condition:
\begin{equation}
\label{Lu's condition}
\sum_i (p_i \cdot u)(h_i \cdot v) = vu
\end{equation}
for every $u,v \in V$.  Then the obvious smash product algebra 
$V \rtimes H'$ is a Hopf algebroid over $V$ with $R$-coring structure and antipode
$\tau$ 
given by ($v \in V, h \in H'$)
\begin{eqnarray}
\tilde{s}(v) & = & v \rtimes 1 \label{one} \\
\tilde{t}(v) & = & \sum_i (p_i \cdot v) \o h_i \label{two} \\
\cop(v \rtimes h) & = & v \rtimes h\1 \o h\2 \label{three} \\
\eps(v \rtimes h) & = & \eps(h)v \label{four} \\
\tau(v \rtimes h) & = & \sum_i (1 \rtimes \tilde{S}(h)) 
\tilde{t}(\tilde{S}^2(h_i) \cdot p_i \cdot v)
\end{eqnarray}
\end{theorem}

\begin{theorem}
The left bialgebroid $T^{\rm op}$ associated to an $H$-Galois extension $A |B$ is
a Hopf algebroid of the type covered in \cite[Theorem 5.1]{Lu}.
\end{theorem}
\begin{proof}
We have seen in example~\ref{ex: tee} that $$T^{\rm op} \cong R \rtimes H^{\rm op}$$
as algebras.  Schauenburg
\cite[3.1]{Sch} computes that the centralizer $R \in \mathcal{YD}^H_H$ 
where $\ract$ denotes the
Miyashita-Ulbrich action of $H$ on $R$,  the coaction $A \to A \o H$
restricts to $R \to R \o H$, and the two intertwine in the following
Yetter-Drinfel'd condition: ($r \in R, h \in H$) 
\begin{equation}
\label{right-right YD condition}
(r \ract h\2)\0 \o h\1 (r \ract h\2)\1 = r\0\! \ract\! h\1 \o  r\1 h\2
\end{equation}  
  
Moreover, the following pre-braided commutativity is satisfied: ($r,r' \in R$)
\begin{equation}
\label{eq: pre-braided commuting}
r' r = r\0 ( r' \ract r\1). 
\end{equation}

Comparing eq.~(\ref{right-right YD condition}) with the left-right Yetter-Drinfel'd
condition \cite[10.6.12]{Mont}, one easily  computes that $$\mathcal{YD}^H_H =
{}_{H^{\rm op}}\!\mathcal{YD}^{H^{\rm op}}$$
when we note that right modules over an algebra correspond exactly to left modules over
its opposite algebra, and that $H^{\rm op}$ has the same coalgebra structure as
$H$ (but with antipode $\overline{S} := \tilde{S}^{-1}$).  
In other words, there are natural actions of $H^{\rm op}$ and its dual on $R$
from the left; the dual acting via the dual of the
coaction (i.e., $p \cdot r = r\0 p(r\1)$) and $H^{\rm op}$ acting
via the Miyashita-Ulbrich action. 
But Majid \cite{Maj} computes that the left-right Yetter-Drinfel'd condition \cite[10.6.12]{Mont}
is equivalent to the anti-commutation relation in the Drinfel'd double
$D(H)= H^{*\, \rm cop} \bowtie H$ (cf.\ \cite{Mont, Kas})  
given by $$(1 \bowtie h)(p \bowtie 1) = (h\1 \rightharpoondown p\2) \bowtie (h\2
\leftharpoondown p\1) $$
where $\rightharpoondown$ and $\leftharpoondown$ denote the right and left coadjoint
actions of $H$ on $H^*$ and $H^*$ on $H$ \cite[10.3.1]{Mont}; whence the left $D(H)$-modules
correspond exactly to left-right Yetter-Drinfel'd modules, or equivalently, 
$$ {}_{D(H^{\rm op})}\!\mathbf{Mod} =  {}_{H^{\rm op}}\!\mathcal{YD}^{H^{\rm op}}.$$

Then $R$ is a left $D(H^{\rm op})$-module; since the coalgebra structure of
$D(H^{\rm op})$ is just $H^* \o H$, we see this action is measuring as well. 
It follows that $R \rtimes H^{\rm op}$ in example~\ref{ex: tee} is a
 smash product $V \rtimes H$ of the type
satisfying the conditions  in Theorem~\ref{thm-Lu}  with
$V = R$ and $H' = H^{\rm op}$, for $D(H^{\rm op}$ has the
R-matrix $\sum_i (1 \bowtie h_i) \o (p_i \bowtie 1)$ (where $\sum_i p_i(x)h_i = x$ for
each $x \in H$ and $p_i(h_j) = \delta_{ij}$), so we compute using Eq. (\ref{eq: pre-braided commuting}):
\begin{eqnarray*}
\sum_i (p_i \cdot y)(h_i \cdot x) & = & y\0 p_i(y\1) (x \ract h_i) \\
& = & y\0 (x \ract y\1) = xy
\end{eqnarray*}

Finally we compute that the bialgebra structure on $R \rtimes H^{\rm op}$
coming from $T^{\rm op}$ in example~\ref{ex: tee} is identical with
that of eqs.\ (\ref{one})-(\ref{four}).  

$$\tilde{t}(r) = \sum_i (p_i \cdot r) \o h_i = \sum_i r\0 \o h_i p_i(r\1) = r\0 \o r\1.$$
This and the other $R$-coring structures are then clearly the same.

We conclude that $T^{\rm op}$ is a Hopf algebroid with antipode $\tau$ on 
$R \rtimes H^{\rm op}$ given by
\begin{equation}
\tau(r \rtimes h) = (1 \rtimes \overline{S}(h))
 (r\0 \ract \overline{S}^2(r\1))\0
  \rtimes (r\0 \ract \overline{S}^2(r\1))\1. \qed
\end{equation}
\renewcommand{\qed}{}\end{proof}


\section{A Galois correspondence for H-separable extensions of simple rings}

  Although Hopf-Galois
extensions in general lack a main theorem of Galois theory \cite{OZ},
we expose results of Sugano in light of obtaining a Galois correspondence
for a depth two cousin of Hopf-Galois extensions, namely H-separable extensions. Their 
definition and part of the proposition below are due to \cite[Hirata]{Hi,H2}.
We will require the Hopf algebroids introduced for H-separable extensions in \cite{K2003}. 
We must eventually narrow our focus to certain H-separable extensions of  \textit{simple rings}, which
 in this section will denote  rings with no proper
two-sided ideals; such a ring is not necessarily artinian or finite dimensional over
a field.   Again let $B$ be a subring of $A$ with centralizer subring  $R$, endomorphism ring 
$S = \End\,_BA_B$ and ring $T = (A \o_B A)^B$.

\begin{lem&def}
$A | B$ is H-separable if $A \o_B A \oplus * \cong \oplus^n A$ as $A$-$A$-bimodules.
Equivalently, $A | B$ is H-separable if there are elements $e_i \in (A \o_B A)^A$
and $r_i \in R$ (a so-called H-separability system) such that
\begin{equation}
1 \o 1 = \sum_i r_i e_i.
\end{equation}
\end{lem&def}

We note that $e_i \in T$, and for $a, a' \in A$
$$ a \o a' = \sum_i e_i \rho_{r_i}(a)a' = \sum_i a\lambda_{r_i}(a')e_i, $$
whence $e_i, \lambda_{r_i}$ is a right D2 quasibasis and $e_i, \rho_{r_i}$ is a
left D2 quasibasis for $A | B$.

For example, given an Azumaya algebra $D | K$ and an arbitrary $K$-algebra $B$
then $A := D \o_K B$ is an H-separable extension of $B$ 
\cite{Hi}. If $B$ is a type $II_1$ factor and $D = M_n(\C)$, this example
covers all H-separable finite Jones index subfactors $B \subseteq A$
 by Proposition~\ref{prop: characterization}(2)
and Proposition~\ref{prop: sugano} below.  

We next let $Z$ denote the center of $A$.  

\begin{prop}
\label{prop: characterization}
If $A | B$ is an H-separable extension, then
\begin{enumerate}
\item $R$ is f.g. projective $Z$-module;
\item $A \o_Z R^{\rm op} \cong \End A_B$ via $a \o r \mapsto \lambda_a \rho_r$;
\item $R \o_{Z} R^{\rm op} \cong S$ via $r \o r' \mapsto \lambda_r \rho_{r'}$
is an isomorphism of bialgebroids;
\item $A \o_B A \cong {\rm Hom}_Z(R, A)$ via $a \o a' \mapsto \lambda_a \rho_{a'}$.
\item $T^{\rm op} \cong \End_Z R$ via $t \mapsto \lambda(t^1) \rho(t^2)$
is an isomorphism of bialgebroids.
\end{enumerate}
Conversely, if $A_B$ is f.g.\ projective,
 the first two conditions imply that $A$ is an H-separable extension
of $B$.  
\end{prop}
\begin{proof}
We offer some short alternative proofs to these facts. $R_Z$ is f.g. projective
since for each $r \in R$, we note that $r = r_i e_i^1 r e_i^2$ 
where summation over $i$ is understood
and for each $i$, $r \mapsto e_i^1 r e_i^2$ defines a map in $\Hom_Z(R,Z)$.   

The inverse $\End A_B \to A \o_{Z} R^{\rm op}$ to the ring homomorphism above is given by 
$f \mapsto  f(e_i^1)e_i^2 \o r_i, $
since $f(e^1_i)e^2_iar_i = f(ae_i^1)e_i^2r_i = f(a)$
 ($a \in A$), while
$ ae_i^1 r e^2_i \o_Z r_i = a \o e_i^1 r e_i^2 r_i = a \o r$ for
$r \in R$ and $a \in A$. 

The ring isomorphism $R^e \cong S$ follows from noting the previous isomorphism is
an $A$-$A$-bimodule morphism. That this ring isomorphism preserves
the $R$-bialgebroid structures on $S$ (cf.\ Section~2)
with respect to Lu's $R$-bialgebroid structure on $R \o_Z R^{\rm op}$ follows
from \cite[5.1]{K2003}.

The inverse ${\rm Hom}_Z (R, A) \to A \o_B A$ to the ring homomorphism above
is given by $g \mapsto g(r_i)e_i$,
since for each $a \o a' \in A \o_B A$, $a r_i a' e_i = a r_i e_i a' = a \o a'$,
while for each $r \in R$, 
$g \in \Hom(R,A)$, $g(r_i) e_i^1 r e_i^2 = g(r_i e_i^1 r e_i^2) = g(r)$.

The ring isomorphism $\End_Z R \cong T^{\rm op}$ follows from noting that the previous
mapping is an $A$-$A$-bimodule homomorphism; this preserves the Lu and depth two
$R$-bialgebroid structures by \cite[5.2]{K2003}.

The converse follows from the general fact that $$H := \Hom (A \o_B A_A, A_A) \cong \End A_B.$$
Since $R$ is f.g.\ projective over $Z$, $\mathcal{E} := \End A_B \cong A \o_Z R$ implies 
that $\mathcal{E}$ is centrally projective over $A$: 
$$\mathcal{E} \oplus * \cong A \oplus \cdots \oplus A,$$
whence the same is true of 
$A \o_B A \cong \Hom ({}_AH, {}_AA)$, which follows from  
$A \o_B A_A$ being finite projective.  
\end{proof}

The $R$-bialgebroid $S$ is in fact an Hopf algebroid since the obvious antipode
on $R^e$ (cf.\ \cite[Lu]{Lu}) is transferred via part (3) of the proposition
\cite{K2003}. 

We prove a lemma relevant to section~2 but independent of the rest of this section.
\begin{lemma}
\label{lemma-seeme}
If $A/B$ is H-separable and $A_B$ is balanced, then $C_A(C_A(B)) = B$.
\end{lemma}
\begin{proof}
Since $R^e \cong \mathcal{S}$, we have for each
$\alpha \in S$ we find $r^1 \o r^2 \in R^e$ such that $\alpha = \lambda(r^1) \rho(r^2)$.
Then for $t \in C_A(R)$:
$$ \alpha \ract t = \lambda(r) \rho(r) t = rr't = \alpha(1) t. $$
So $t \in A^S$, an invariant under the action, whence $t \in B$ since $A_B$ is balanced
\cite[Theorem 4.1]{KSz}.
\end{proof} 

The proposition and theorem below are due to Sugano, recapitulated below
in a hopefully useful expository manner. 

\begin{prop}[Sugano \cite{Su82}]
\label{prop: sugano}
Suppose $B$ is a simple ring and subring of $A$.  Then $A$ is a right f.g.\ projective H-separable extension
of $B$ if and only if
\begin{enumerate}
\item $A$ is a simple ring,
\item $C_A(C_A(B)) = B$, and
\item $C_A(B)$ is a simple finite dimensional $Z$-algebra.  
\end{enumerate}
\end{prop}
\begin{proof} 
($\Rightarrow$) Since $\Hom (A_B, B_B) \neq 0$ and $B$ has no non-trivial ideals, 
the trace ideal for $A_B$ is $B$, so $A_B$ is a generator. Let
$f_i \in \Hom (A_B,B_B)$, $a_i \in A$ such that $\sum_i f_i(a_i) = 1_A$.  Then the
inclusion $\iota: B \to A$ is split as right $B$-module mapping by $a \mapsto
\sum_i f_i(a_ia)$.  Let $e: A_B \to B_B$ be a projection.  Given a two-sided
 ideal $I \subset A$,
we have $$I = (I \cap B)A$$ since $x = \sum_i e(r_i xe_i^1) e_i^2$ for $x \in I$ and
H-separability system $\{ e_i, r_i \}$; but $e(r_i x e_i^1) \in I \cap B$ by
Proposition~\ref{prop: characterization}(2). Then $B$ simple implies $I = 0$,
whence $A$ is simple. 

Clearly, $B \subseteq C_A(R)$ where $R = C_A(B)$.  Let $v \in C_A(R)$
and $\phi$ denote the isomorphism in Proposition~\ref{prop: characterization}(4).
Then $\phi(v \o 1)(r) = \phi(1 \o v)(r)$ for every $r \in R$, whence
$1 \o v = v \o 1$ in $A \o_B A$.  Applying the projection $e$, we arrive
at $v = e(v) \in B$, whence $B = C_A(R)$.  

Since $A_B$ is a progenerator, $\End A_B$ is also a simple ring by Morita theorems.
Then $A \o_Z R$ is simple. Since $Z$ is a field by Schur's lemma, it follows that
$R$ is a simple (finite dimensional) $Z$-algebra.

($\Leftarrow$)  The map in  Proposition~\ref{prop: characterization}(2), call
it $\psi$, always
exists although it may not be an isomorphism. By conditions (1) and (3) however,
$\psi$ is a monomorphism from $\Lambda := A \o_Z R$ into $\End A_B$. It suffices to
 show that $\psi$ is an isomorphism and $A_B$ is f.g.\ projective by the converse
in Proposition~\ref{prop: characterization}.   If $C$ is the center of $R$,
it follows from $C \o_Z C$ being a Kasch ring and $R$ being a $C$-separable
algebra that $\Hom (A_{\Lambda}, \Lambda_{\Lambda})
\neq 0$ \cite{Su82}, whence $A_{\Lambda}$ is a left generator, therefore right $\End {}_{\Lambda}A$-f.g.\
projective.  But $\End {}_{\Lambda}A = {\rm Hom}_{A-R}\, (A,A) \cong C_A(R) = B$,
so $A_B$ is f.g.\ projective. Again, $B$ is simple and
$\Hom (A_B, B_B) \neq 0$ implies that $A_B$ is also a generator.
It follows from Morita theorems that $\Lambda \cong \End A_B$
via $\psi$.   
\end{proof}  

In \cite{K2001} a (right) HS-separable extension $A | B$ is defined to be H-separable
such that the natural module $A_B$ is a progenerator.  What we have then seen
above is that a right f.g.\ projective
 H-separable extension $A | B$, where $B$ is simple, is HS-separable. 

\begin{theorem}[Sugano \cite{Su82,Su87}]
Suppose $A$ is an HS-separable extension of a simple ring $B$.  Then the class
of simple $Z$-subalgebras $V$ of the centralizer $R$ is in one-to-one correspondence
with the class of intermediate simple subrings $B \subseteq D \subseteq A$
where $A_D$ is f.g.\ projective, via the centralizer in $A$:
  $D \mapsto C_A(D)$ with inverse $V \mapsto C_A(V)$. Moreover, $A$ over each
such intermediate simple ring $D$ is an HS-separable extension. 
\end{theorem}
\begin{proof}
Of course, $A$ is a simple ring by proposition.   
Given $D$ as in the theorem, we show $D$ is a right relatively separable extension
of $B$ in $A$, i.e., the multiplication mapping $\mu: A \o_B D \to A$ is split
as an $A$-$D$-bimodule epi. For then $$ A \o_D A \oplus * \cong A \o_B D \o_D A  =
A \o_B A, $$
as $A$-$A$-bimodules, the latter being isomorphic itself to a direct summand of
$A \oplus \cdots \oplus A$; whence $A | D$ is H-separable, in fact HS-separable
since $D$ is simple. It follows from the proposition then that $C_A(D)$ is a
simple $Z$-algebra with $C_A(C_A(D)) = D$, which yields half of the theorem.

To show that $D$ is right relatively separable extension of $B$ in $A$,
\cite{Su87} shows by other means from the hypotheses that $A | D$ is H-separable,
hence $C_A(D)$ is simple: not surprisingly then, $R$ is a Frobenius extension of $C_A(D)$,
so $D$ is a
Frobenius extension of $B$ via an isomorphism, say $\eta$, given in Proposition~\ref{prop: characterization}(4)
\cite[Theorem 3]{Su87}. Let $\{ E: D \to B, x_i, y_i \}$ be a Frobenius system for
$D | B$.  Consider the two-sided ideal $I := \sum_i x_i R y_i$ in $C_A(D)$.  If $I = 0$,
then by Proposition~\ref{prop: characterization}(4)
 $\sum_i x_i \o y_i = 0$ in $A \o_B A$, whence in $D \o_B D$ since $A_B$ and ${}_BD$ are
 flat modules.  
It follows that $\sum_i E(x_i)y_i = 0$, which contradicts $\sum_i E(dx_i)y_i = d$
for all $d \in D$.  Then $ = \sum_i x_i R y_i = C_A(D)$.  It follows that
there is $r \in R$ such that $\sum_i x_i r \o y_i \in A \o_B D$ is a right
relative separability element which yields a splitting for $\mu$. 

The other half of the theorem depends on showing that $C_A(V) = D'$ is a simple
ring, which clearly is intermediate to subring $B$ and over-ring $A$, and furthermore
$C_A(D') = V$ as well as $A_{D'}$ being f.g.\ projective.  Since $A$ is a right $B$-generator,
$A$ is left-f.g.\ projective over $\End A_B \cong A \o_Z R$.  But $R$ is a f.g.\
projective $V$-module, whence $A$ is f.g.\ projective left $\Omega := A \o_Z V$-module.
Since $\Omega$ is a simple ring, ${}_{\Omega}A$ is also a generator, so $A_{D'}$ is
a progenerator module and $D'$ is a simple ring by Morita theorems, since $\End {}_{\Omega}A \cong C_A(V)$ via
$f \mapsto f(1)$. 

Now let $V' := C_A(D') \supseteq V$.  Clearly there is a mapping $A \o_Z V' \to \End A_{D'}$
as in Proposition~\ref{prop: characterization}(2), which forms commutative squares
with two other such mappings $A \o_Z V \stackrel{\cong}{\to} \End A_{D'}$
and $A \o_Z R \stackrel{\cong}{\to} \End A_B$.  These squares are joined by inclusions,
which forces $\psi$ to be an isomorphism
 and $A \o V = A \o V'$ over the
field $Z$, whence $V = V'$.  Then $A | D'$ is HS-separable
and the correspondence in the theorem is one-to-one.  
\end{proof} 
 
We are now in a position to establish a Galois correspondence between intermediate
simple rings of $A | B$ and Hopf subalgebroids of $S$ over simple subalgebras
of $R$.  The one-to-one correspondence below bears a resemblance to the Jacobson-Bourbaki correspondence
for division rings.  

\begin{theorem}
Given an HS-separable extension of simple rings $A | B$, there is a one-to-one
correspondence between  intermediate simple rings $D$ such that $A_D$ is f.g.\
projective and Hopf subalgebroids $H$ of $S$ minimal over simple subalgebras $V \subseteq R$.  
The Galois correspondence is given by $D \mapsto \End {}_DA_D$, a Hopf algebroid
over $C_A(D)$, with inverse given by $H \mapsto A^H$, the fixed points under the
canonical action of $H$.  
\end{theorem}
\begin{proof} 
Given a simple intermediate ring $D$ such that $A_D$ is finite projective,
we have seen that $A | D$ is an HS-separable extension, hence a depth two right
balanced extension of the type considered in \cite[Section~4]{KSz}.  It follows
that $H := \End {}_DA_D$ is a left bialgebroid over $C_A(D)$ such
that $A^H = D$ under the action given in Eq.~\ref{eq: alpha-action};
with antipode and Lu Hopf algebroid structure \cite{K2003} from
Proposition~\ref{prop: characterization}(3). There are clearly  inclusions
$H \subseteq S$ and another $C_A(D) \subseteq R$ which together show $H$ to be 
a Hopf subalgebroid of $S$ minimal over $C_A(D)$. 
Of course, $C_A(D)$ is a simple $Z$-algebra
and $C_A(C_A(D)) = D$ by
Proposition~\ref{prop: sugano}.  

Conversely, given a Hopf subalgebroid $H$ of $S$ minimal over
a simple subalgebra $V \subseteq R$, we let $D' = C_A(V)$,
an intermediate ring between $B$ and $A$ which is simple with $A | D'$
an HS-separable extension by the last proposition.  Now under identification
of $S$ with $R \o_Z R^{\rm op}$, we note that 
$V \o_Z V^{\rm op} \subseteq H$ since $s(v) = v \o 1$ and $t(v') = 1 \o v'$
for $v,v' \in V$, while $s(v)t(v') = v \o v' \in H$ as well. 
Since $V \o V^{\rm op}$ is a $V$-bialgebroid and obviously a subbialgebroid of $H$,
the minimality condition forces $H = V \o V^{\rm op}$.  
Since $\End {}_{D'}A_{D'} = V \o V^{\rm op} = H$ and $A | D'$ is 
depth two right balanced, it follows from \cite{KSz} that $A^H = D'$.  
Therefore the correspondence in the theorem is one-to-one.    
\end{proof}

\end{document}